# Parasitic Numbers at Arbitrary Base

**Anatoly A. Grinberg**

**anatoly_gr@yahoo.com**

## ABSTRACT

**A natural number is called λ-parasitic number if it is multiplied by integer λ as the rightmost digit moves to the front. Full set of these numbers are known in the decimal system. Here, a formula to analytically generate parasitic numbers in any basis was derived and demonstrated for the number systems t = 3, 4, 5, 8, 10 and 16. It allows to generate parasitic numbers with given numbers of periods. The formula was derived on the basis of parasitic number definition, but not on their cyclic property.**

•

## I. Introduction.

There are integer numbers, which increase by times **λ** by displacement its lower rank digit $b_0$ to the highest one. In mathematical literature, they are known by not very flattering title as the parasitic numbers (PN). They attract public attention under Dyson's number name when they have been proposed to the public as a challenge to find them **[1].** Full set of these numbers are well known in decimal number system. Some number in duodecimal system also are published **[2].**

Due to nature of this numbers the cyclical number are created by their geometrical «docking». Therefore the cyclical properties of rational fractions is used to create PN. Here we derived general formula using direct definition of the PN, but not their cyclical nature. This gives us a unified method to generate parasitic numbers with given number of the periods in any number system. One period constitute minimal PN. The number of PN digits **N** is the only unknown parameter in the formula. The number **N** has to be found from congruence. The module of the congruence is determined by base number **t** and multiplicative parameter **λ.** We presented parasitic numbers in explicit form for number bases **t** equal **to 3, 4, 5, 8,10,16 .** We also constructed the tables that show parameter **N** as function of digit $\boldsymbol{b_0}$ and multiplicative factor **λ,** for values of **t** in the range **[3 – 20]**. Using these tables a calculation of any **PN** is a straight forward procedure.

## II. The formulation of the equation.

**N**-digit number $\boldsymbol{P_\lambda^N}$ in the t-base system has form

$$\boldsymbol{P_\lambda^N = b_{N-1}\, t^{N-1} + b_{N-2}\, t^{N-2} + \cdots \ldots . b_1 t + b_0 = Bt + b_0} \qquad \textbf{(II.1)}$$

where digits $b_k$ satisfy inequality $|\boldsymbol{b_k}| \le (\boldsymbol{t-1})$.

By moving the digit $\boldsymbol{b_0}$ to the front (i.e. leftmost) place we obtain

$$\boldsymbol{Q_\lambda^N = b_0\, t^{N-1} + b_{N-1}\, t^{N-2} + b_{N-2}\, t^{N-3} + \cdots b_1 = b_0 t^{N-1} + B} \qquad \textbf{(II.2)}$$

In agreement with the definition of the problem $\boldsymbol{Q_\lambda^N}$ is $\boldsymbol{\lambda}$-times larger than $\boldsymbol{P_\lambda^N}$. Thus, we have equation

$$\boldsymbol{Q_\lambda^N = \lambda\, P_\lambda^N} \qquad \textit{\textbf{(II.3)}}$$

from which it follows that

$$\boldsymbol{b_0 t^{N-1} + B = \lambda(Bt + b_0)} \qquad \textbf{(II.4)}$$

Substituting **B** from **(II.4)** into **(II.1)** and **(II.2)** we find

$$\boldsymbol{P_\lambda^N(\mathbf{b_0}) = b_0\, \frac{(t^N - 1)}{(\lambda\, t - 1)}} \qquad \textbf{(II.5)}$$

$$\boldsymbol{Q_\lambda^N(\mathbf{b_0}) = \lambda\, b_0\, \frac{(t^N - 1)}{(\lambda\, t - 1)}} \qquad \textbf{(II.6)}$$

Parameter **N** is an unknown in **(II.5), (II.6).** Because $P_\lambda^N$ and $\boldsymbol{Q_\lambda^N}$ must be integers this is the only condition from which **N** can be found. The divisibility of numerator of **(II.5)** by $\boldsymbol{p} = (\boldsymbol{\lambda\, t - 1})$ can be expressed in the form of congruence

$$\mathbf{b_0}\, (\boldsymbol{t^N} - \mathbf{1}) \equiv \mathbf{0} \qquad (\mathbf{mod}\, (\boldsymbol{\lambda}\, \mathbf{t} - \mathbf{1})) \qquad \textbf{(II.7)}$$

In **Appendix I** it is shown, that congruence **(II.7)** can be transformed into following congruence

$$\mathbf{b_0}\, (\boldsymbol{\lambda^N} - \mathbf{1}) \equiv \mathbf{0} \qquad (\mathbf{mod}\, (\boldsymbol{\lambda}\, \mathbf{t} - \mathbf{1})) \qquad \textbf{(II.8)}$$

The latter was deduced using General Divisibility Criteria (**GDC**) presented in **ref.[3].** The two congruences are equivalent in the considered here problem and can be used interchangeable. It is usually better to use the congruence with smaller base power. It worth to noting that, if the module $\boldsymbol{p} = (\boldsymbol{\lambda}\, \mathbf{t} - \mathbf{1})$ is co-prime to $\mathbf{b_0}$ **,** then the latter can be removed from congruence. In the opposite case, $\mathbf{b_0}$ can also be removed, but the module must be divided by the greatest common divisor of $\mathbf{b_0}$ and **p**. Therefore the classical theory of power residue elaborated by Gauss [4] can be utilized. It is important that the power bases **t** and parameter **λ** are co-prime to congruence module **p**.

Because all parameter we used in **(II.6)** are expressed in decimal base, the calculated number $P_\lambda^N$ also expressed in this base. They have to be transformed into base **t**. Only in this system will they have N digits and the properties of the PN numbers.

Factor **λ** cannot be larger than **(t-1)** because the length of the numbers $P_\lambda^N$ and $Q_\lambda^N$ must be the same. It is evident that at given value of **λ,** formula **(II.5)** generates **(t-1)** numbers of $P_\lambda^N$ in accordance with possible values of digit $b_0$**.** However, the numbers of digits in $P_\lambda^N$ and $Q_\lambda^N$ will not be the same for all $b_0$. It will be the same only if inequality $b_0 \geq \lambda$ is satisfied**.** Solutions with unequal number of digits must be removed as not satisfying the condition of the problem. The number of solutions at any given value of $\lambda$ decreases with $\lambda$ increase. In particularly, at $\lambda = t-1$ only one parasitic number $P_{(t-1)}^N(t-1)$ exists**.**

The $P_\lambda^N$ and $Q_\lambda^N$ are determined by three parameters: $t$, $b_0$ *and* $\lambda$**. N** is found from appropriate congruence.

At the chapter **IX** we presented the tables of **N** as function of $b_0$ *and* $\lambda$ for bases t in the interval of **t = [3 - 20].** The tables were constructed by solution of congruence **(II.6**), **(II.7).** It is straight forward to find $P_\lambda^N$ and $Q_\lambda^N$ using these tables and **Eqs.(II.5)** and **(II.6).**

The explicit form of the numbers $P_\lambda^N$, $Q_\lambda^N$ are presented below for different values of digit $b_0$ , factor **λ** and base **t**. The number are shown in decimal and appropriate t-base systems. The number for non-decimal base were written in brackets with corresponding t subscript.

Because $P_\lambda^N$ and $Q_\lambda^N$ are multiple of $b_0$ they are presented in the form $P_\lambda^N(b_0) = b_0 * P_\lambda^N(1)$ with $P_\lambda^N(1)$ indicated in decimal and t-base representation The number $P_\lambda^N$ *and* $Q_\lambda^N$ are written in explicit form when only one solution exists at particular $\lambda$ and $N$.

The results of these calculations are presented in the following schema. On the first line we show: the value of the multiplicative parameter **λ**, the range of the possible $b_0$ values, appropriate congruence, and following from this the number **N.** On the following lines, the numbers $P_\lambda^N$ and $Q_\lambda^N$ are presented at decimal and t-base formats.

## III. Ternary Number System (t=3)

In ternary number system parameter λ has only one value, λ = 2 and only one number $P_2^4(2)$ exists**:**

$$\lambda = 2 \ | \ b_0 = [2] \ | \ [3^4 \equiv 1 \ (mod\ 5)] \ | \ N = 4$$

$$P_2^4(2) = 32 = (1012)_3 ; \quad Q_2^4(2) = 64 = (2101)_3 \qquad \textbf{(III.1)}$$

## IV. Quaternary Number System (t=4)

$$\lambda = 2 \ | \ b_0 = [2-3] \ | \ [4^3 \equiv 1 \ (\mathrm{mod}\ 7)] \ | \ N = 3$$

$$P_2^3(b_0) = 9b_0 = (21)_4 \, b_0$$

$$P_2^3(2) = 18 = (102)_4\,; \quad Q_2^3(2) = 36 = (210)_4$$

$$P_2^3(3) = 27 = (123)_4\,; \quad Q_2^3(3) = 54 = (312)_4 \qquad \textbf{(IV.1)}$$

$$\lambda = 3 \;|\; b_0 = [3] \;|\; \left[4^5 \equiv 1 \;(mod\;11)\right] \;|\; N = 5$$

$$P_3^5(3) = 279 = (10113)_4\,; \quad Q_3^5(3) = 83 = (31011)_4 \qquad \textbf{(IV.2)}$$

## *V.* Quinary Number System ( t=5)

$$\lambda = 2 \;|\; b_0 = [2, 4] \;|\; [5^6 \equiv 1\;(mod\;9)] \;|\; \mathrm{N} = 6$$

$$P_2^6(b_0) = 1736 b_0 = (23421)_5\, \mathrm{b}_0$$

$$P_2^6(2) = 3472 = (102342)_5\;;\; Q_2^6(2) = 6944 = (210234)_5$$

$$P_2^6(4) = 6944 = (210234)_5\,;\; Q_2^6(4) = 13888 = (421023)_5$$

$$\lambda = 2 \;|\; b_0 = [3] \;|\; [5^2 \equiv 1\;(mod\;3)] \;|\; \mathrm{N} = 2$$

$$P_2^2(3) = 8 = (13)_5\,; \quad Q_2^2(3) = 16 = (31)_5 \qquad \textbf{(IV.3)}$$

$$\lambda = 3 \;|\; b_0 = [3 - 4] \;|\; [5^6 \equiv 1\;(mod\;14)] \;|\; \mathrm{N} = 6$$

$$P_3^6(b_0) = 1116\, b_0 = (1, 3, 4, 3, 1)_5\, b_0$$

$$P_3^6(3) = 3348 = (101343)_5\;; Q_3^6(3) = 10044 = (310134)_5$$

$$P_3^6(4) = 4464 = (120324)_5\;; Q_3^6(4) = 13392 = (412032)_5 \qquad \textbf{(IV.4)}$$

$$\lambda = 4 \;|\; b_0 = [4] \;|\; [5^9 \equiv 1\;(mod\;19)] \;|\; \mathrm{N} = 9$$

$$P_4^9(4) = \quad 411184 = (101124214)_5$$

$$Q_4^9(4) = 1644736 = (410112421)_5 \qquad \textbf{(IV.5)}$$

## *VI.* Octal Number System (t=8)

$$\lambda = 2 \mid b_0 = [2-4; 6-7] \mid [\,8^4 \equiv 1 \ (mod\ 15\,)] \mid N = 4$$

$$P_2^4(b_0) = 273\ b_0 = \ (421)_8\, b_0\ ;\ Q_2^4(b_0) = 546\ b_0 = (1042)_8\ b_0$$

$$P_2^4\,(2) = 546 = (1042)_8\,; \quad Q_2^4(2) = 1092 = (2104)_8 \qquad \textbf{(VI.1)}$$

$$\lambda = 2 \mid b_0 = [5] \mid [\,8^2 \equiv 1 \ (mod\ 3\,)] \mid N = 2$$

$$P_2^2(5) = 21 = (25)_8\,; \quad Q_2^2(3) = 42 = (52)_8 \qquad \textbf{(VI.2)}$$

$$\lambda = 3 \mid b_0 = [3-7] \mid [8^{11} \equiv 1 (mod\ 23)] \mid N = 11$$

$$P_3^{11}(b_0) = \ 373475417 b_0 = (2,6,2,0,5,4,4,1,3,1)_8\, b_0$$

$$P_3^{11}(3) = 1\,120,426251 = (10,262,054413)_8$$

$$Q_3^{11}(3) = 3\,361,278753 = (31,026,205441)_8 \qquad \textbf{(VI.3)}$$

$$\lambda = 4\,; b_0 = [4-7];\ \left[8^5 \equiv 1 (mod\ 31)\right]; N = 5$$

$$P_4^5(b_0) = 1057\ b_0 = (2041)_8\ b_0$$

$$Q_4^5(b_0) = 4228\ b_0 = (10204)_8\,\mathrm{b_0} \qquad \textbf{(VI.4)}$$

$$\lambda = 5\,; b_0 = [5-7];\ [8^4 \equiv 1 (mod\ 39)]\,; N = 4$$

$P_5^4(b_0) = 105\, b_0 = (151)_8\, ;\; Q_5^4\,(b_0) = 525\, b_0 = (1015)_8\, b_0$ **(VI.5)**

$\lambda = 6\, ; b_0 = [6-7];\; [8^{23} \equiv 1(mod\, 47)]\, ; N = 23$

$P_6^{23}(b_0) = 2559485326780971313\, b_0 =$

$(127114202562304053446 1)_8\, b_0$

$Q_6^{23}(b_0) = 7535691196068582787 8\, b_0 =$

$(101271142025623040534 46)_8\, b_0$ **(VI.6)**

$\lambda = 7\, ;\; b_0 = [7];\;\; [8^{20} \equiv 1(mod\, 55)]\, ;\;\; N = 20$

$P_6^{23}(7) = 146735464222689615 =$

$(10112362022474404517)_8$

$Q_6^{23}(\,7\,) = 1027148249558827305 =$

$(71011236202247440451)_8$ **(VI.7)**

## *VII.* Decimal Number System (t=10)

$\lambda = 2\;\; | \; b_0 = [2-9]\;\; | \; [10^{18}\; \equiv\; 1\;\; (\;\; mod\; 19\;)]\;\; | \; N = \;\; 18$

$P_2^{18}(b_0) = 52631, 578947, 368421\, b_0$ **(VII.1)**

$\lambda = 3\;\; |\;\;\; b_0 = [3-9]\;\;\; |\;\;\; [10^{28} \equiv 1\;\; (\; mod\; 29\;)]\;\; | \; N = \;\; 28$

$P_3^{28}(b_0) = 344, 827586, 206896, 551724, 137931 b_0$ **(VII.2)**

$\lambda = 4\;\; |\;\;\; b_0 = [4-9]\;\;\; |\;\; [10^6 \equiv 1\;\; (\; mod\; 39\;)]\;\;\; |\;\; N = \;\; 6$

$P_4^6(b_0) = 25641 b_0$ (VII.3)

$\lambda = 5 \mid b_0 = [5, 6, 8, 9] \mid [10^{42} \equiv 1 \ (mod\ 49)] \mid N = 42$

$P_5^{42}(b_0) = 20408, 163265, 306122, 448979, 591836,$

$734693, 877551\ b_0$ (VII.4)

$\lambda = 5 \mid b_0 = [7] \mid [10^6 \equiv 1 \ (mod\ 7)] \mid \mathrm{N} = 6$

$P_5^6(7) = 142857$ ; $Q_5^6(7) = 714285$ (VII.5)

$\lambda = 6 \mid b_0 = [6 - 9] \mid [10^{58} \equiv 1 \ (mod\ 59)] \mid N = 58$

$P_6^{58}(b_0) = 169, 491525, 423728, 813559, 322033,$

$898305, 084745, 762711, 864406, 779661 b_0$ (VII.6)

$\lambda = 7 \mid b_0 = [7 - 9] \mid [10^{22} \equiv 1 \ (mod\ 69)] \mid N = 22$

$P_7^{22}(\mathrm{b}_0) = 144, 927536, 231884, 057971\ b_0$ (VII.7)

$\lambda = 8 \mid b_0 = [8 - 9] \mid [10^{13} \equiv 1 \ (mod\ 79)] \mid \mathrm{N} = 13$

$P_8^{13}(\mathrm{b}_0) = 126582, 278481\ \mathrm{b}_0$ (VII.8)

$\lambda = 9 \mid b_0 = [9] \mid [10^{44} \equiv 1 \ (mod\ 89)] \mid N = 44$

$P_9^{44}(\mathrm{b}_0) = 1, 123595, 505617, 977528, 089887,$

$640449, 438202, 247191\ \mathrm{b}_0$ (VII.9)

**VIII. Hexadecimal Number System** (t=16)

$\lambda = 2;\ b_0 = [2-15];\ \ \left[16^5 \equiv 1\ (mod\ 31)\right]\ ;\ N = 5$

$P_2^5(b_0) = 33825b_0\ =\ (8,4,2,1)_{16}\, b_0$

$P_2^5(2) = 67650 = (1,0,8,4,2)_{16};\ Q_2^5(2) = (2,1,0,8,4)_{16}$ **(VIII.1)**

$\lambda = 3 \mid\ b_0 = [3-15]\ \mid [16^{23} \equiv 1\ (mod\ 47)] \mid N = 23$

$P_3^{23}(b_0) = 1580{,}348986{,}321762{,}053062{,}711775\, b_0 =$

$(5,1,11,3,11,14,10,3,6,7,7,13,4,6,12,14,15,$

$10,8,13,9,13,15)_{16}\, b_0$ **(VIII.2)**

$\lambda = 4\ \mid b_0 = [4-15] \mid\ [16^3 \equiv 1\ (mod\ 63)] \mid N = 3$

$P_4^3(b_0) =\ 975b_0 = (3,12,15)_{16}\, b_0$ **(VIII.3)**

$\lambda = 5\ \mid\ b_0 = [5-15] \mid\ [16^{39} \equiv 1\ (mod\ 79)] \mid N = 39$

$P_5^{39}(b_0) = 1156{,}251295{,}356726{,}992249{,}750658{,}794527,$

$705423{,}997265\, b_0\ =$

$(3,3,13,9,1,13,2,10,2,0,6,7,11,2,3,10,5,4,4,0,12,15,6,$

$4,7,4,10,8,8,1,9,14,12,8,14,9,5,1)_{16}\, b_0$ **(VIII.4)**

$\lambda = 6\ |b_0 = [6-15] \mid [16^9 \equiv 1\ (mod\ 19)] \mid N = 9$

$P_6^9(b_0) = 723362913\, b_0 = (2,11,1,13,10,4,6,1)_{16}\, b_0$ **(VIII.5)**

$\lambda = 7 \mid \; b_0 = [7 - 15] \;\mid\; [\,16^9 \equiv 1\ (mod\ 37\,)] \;\mid\; N = 9$

$$P_7^9(b_0) = 619094385\ b_0\ = (2,4,14,6,10,1,7,1)_{16}\ b_0 \qquad \textbf{(VIII.6)}$$

$\lambda = 8 \mid b_0 = [8 - 15] \mid [\,16^7 \equiv 1\ (mod\ 127\,)] \;\mid\; N = 7$

$$P_5^{39}(b_0) = 2113665\ b_0 = (1,14,3,12,7,8,15)_{16}\ b_0 \qquad \textbf{(VIII.7)}$$

$\lambda = 9 \mid b_0 = [9,10.12,14,15] \mid \left[\,16^{15} \equiv 1\ (mod\ 143\,)\right] \mid N = 15$

$P_9^{15}(b_0) =\ 8062388144103825\ b_0 =$

$$(1,12,10,4,11,3,0,5,5,14,14,1,9,1)_{16}\ b_0 \qquad \textbf{(VIII.8)}$$

$\lambda = 9 \;\mid\; b_0 = [11] \mid [\,16^3 \equiv 1\ (mod\ 13\,)] \mid N = 3$

$$P_9^3(11) = 315 = (1,3,11)_{16}\,;\ Q_9^3(11) = 2835 = (11,1,3)_{16} \qquad \textbf{(VIII.9)}$$

$\lambda = 9 \mid b_0 = [13] \;\mid\; \left[\,16^5 \equiv 1\ (mod\ 11\,)\right] \mid N = 5$

$P_9^5(13) = 95325\ \ = (1,7,4,5,13)_{16};$

$$Q_9^5(13) = 857925 = (13,1,7,4,5)_{16} \qquad \textbf{(VIII.10)}$$

$\lambda = 10 \;\mid b_0 = [10 - 15] \mid [\,16^{13} \equiv 1\ (mod\ 159\,)] \mid N = 13$

$P_{10}^{13}(b_0) =\ 28324525958305\ b_0 =$

$$(1,9,12,2,13,1,4,14,14,4,10,1)_{16}\ b_0 \qquad \textbf{(VIII.11)}$$

$\lambda = 11 |\ b_0 = [11 - 13] \mid \left[\ 16^{15} \equiv 1\ (mod\ 175\ )\right] \mid N = 15$

$P_{11}^{15}(b_0) = \ 6588122883467697\ b_0 =$

$$(1,7,6,7,13,12,14,4,3,4,10,9,11,1)_{16}\ b_0 \qquad \text{(VIII.12)}$$

$\lambda = 11\ \mid b_0 = [14]\ \left|\left[\ 16^5 \equiv 1\ (mod\ 25\ )\right]\ \right|\ N = 5$

$P_{11}^{5}(14) = 83886 \quad = (1,4,7,10,14)_{16}$

$$Q_{11}^{5}(14) = 922746 = (14,1,4,7,10)_{16} \qquad \text{(VIII.13)}$$

$\lambda = 11 \mid b_0 = [15] \mid\ [\ 16^3 \equiv 1\ (mod\ 35\ )]\ \mid N = 3$

$P_{11}^{3}(15) = 351 = (1,5,15)_{16}$

$$Q_{11}^{3}\ (15) = 3861 = (15,1,5)_{16} \qquad \text{(VIII.14)}$$

$\lambda = 12 \mid\ b_0 = [12 - 15] \mid \left[\ 16^{95} \equiv 1\ (mod\ 191\ )\right] \mid\ N = 95$

$P_{12}^{95}(b_0) = 12893, 326634, 945837, 438572, 984325,$

$963224, 412656, 982745, 571154, 014380, 986061, 598730,$

$946170, 553210, 541317, 491781, 713233, 248242, 143425\ b_0 =$

$(1,5,7,1,14,13,3,12,5,0,6,11,3,9,10,2,2,13,9,2,1,8,2,0,2,$

$10,14,3,13,10,7,8,10,0,13,6,7,\ 4,3,0,4,0,5,5,12,7,11,4,$
$15,1,4,1,10,12,14,6,8,8,11,6,4,8,6,0,8,0,10,11,8,15,6,$
$$9,14,2,8,3,5,9,12,13,1,1,6,12,9,0,\ 12,1)_{16}\ b_0 \qquad \text{(VIII.15)}$$

$\lambda = 13 \mid b_0 = [13 - 14] \mid [\ 16^{33} \equiv 1\ (mod\ 207\ )] \mid N = 33$

$P_{13}^{33}(b_0) = 2630201869920297302132364115 4146$

$335185\, b_0 = (1,3,12,9,9,5,10,4,7,11,10,11,14,7,$

$4,4,0,4,15,2,6,5,6,9,1,14,14,10,15,9,13,1\,)_{16}\, b_0$ **(VIII.16)**

$\lambda = 13 \mid b_0 = [15] \mid [\,16^{11} \equiv 1\ (mod\ 69\,)] \mid N = 11$

$P_{13}^{11}(15) = 1,274796,090175 =$

$(1,2,8,12,15,12,4,10,3,3,15)_{16}$

$Q_{13}^{11}(15) = 16572349172275 =$

$(15,1,2,8,12,15,12,4,10,3,3)_{16}$ **(VIII.17)**

$\lambda = 14 \mid b_0 = [14-15] \mid [\,16^{37} \equiv 1\ (mod\ 223)] \mid N = 37$

$P_{14}^{37}(b_0\,) = 1,600053,467159,147848,720051,535257,$

$281543,029985\, b_0 = (1,2,5,14,2,2,7,0,8,0,9,2,15,1,1,$

$3,8,4,0,4,9,7,8,8,9,12,2,0,2,4,11,12,4,4,14,1)\,_{16}\, \mathrm{b_0}$ **(VIII.18)**

$\lambda = 15 \mid b_0 = [15] \mid [16^{119} \equiv 1\ (mod\ 239\,)] \mid N = 119$

$P_{15}^{119}(15) = 12245,352577,076329,359498,943108,$

$498298,793469,919492,671173,447407,780960,$

$873581,611479,206145,869401,554903,864949,$

$633212,334419,743552,152646,348114,363043,$

$826975 = (1,0,1,1,2,3,5,8,14,7,5,13,3,0,3,3,6,10,0,10,$

$11,6,1,7,9,0,9,10,3,14,2,0,2,2,4,6,11,1,12,14,11,10,6,$

$0,6,6,13,4,1,5,6,12,2,15,2,1,3,4,7,12,4,0,4,4,8,13,6,3,$

$9,13,7,4,12,0,12,13,10,8,2,10,13,8,5,14,4,2,6,8,15,8,8,$

$$0, 8, 9, 1, 10, 12, 7, 3, 10, 14, 9, 8, 1, 9, 11, 5, 0, 5, 5, 11, 0, 11, 12,$$
$$8, 4, 13, 1, 15)_{16} \quad \textbf{(VIII.19)}$$

## IX. Tables representing parameter N of PN as function of transposed digit $b_0$ and multiplicative factor λ for values of t = [3-20].

| | t=3 | |
|---|---|---|
| $b_0$ | | |
| **2** | 4 | |
| | **2** | λ |

| | t=4 | | |
|---|---|---|---|
| $b_0$ | | | |
| **3** | 3 | 5 | |
| **2** | 3 | | |
| | **2** | **3** | λ |

| | | t=5 | | |
|---|---|---|---|---|
| $b_0$ | | | | |
| **4** | 6 | 6 | 9 | |
| **3** | 2 | 6 | | |
| **2** | 6 | | | |
| | **2** | **3** | **4** | λ |

| | | | | t=6 | |
|---|---|---|---|---|---|
| $b_0$ | | | | | |
| **5** | 10 | 16 | 11 | 14 | |
| **4** | 10 | 16 | 11 | | |
| **3** | 10 | 16 | | | |
| **2** | 10 | | | | |
| | **2** | **3** | **4** | **5** | λ |

| | | | t=7 | | | |
|---|---|---|---|---|---|---|
| $b_0$ | | | | | | |
| **6** | 12 | 4 | 3 | 16 | 40 | |
| **5** | 12 | 2 | 9 | 16 | | |
| **4** | 12 | 4 | 9 | | | |
| **3** | 12 | 4 | | | | |
| **2** | 12 | | | | | |
| | **2** | **3** | **4** | **5** | **6** | λ |

| | | | | t=8 | | | |
|---|---|---|---|---|---|---|---|
| $b_0$ | | | | | | | |
| **7** | 4 | 11 | 5 | 4 | 23 | 20 | |
| **6** | 4 | 11 | 5 | 4 | 23 | | |
| **5** | 2 | 11 | 5 | 4 | | | |
| **4** | 4 | 11 | 5 | | | | |
| **3** | 4 | 11 | | | | | |
| **2** | 4 | | | | | | |
| | **2** | **3** | **4** | **5** | **6** | **7** | λ |

| | | | | | t=9 | | | |
|---|---|---|---|---|---|---|---|---|
| $b_0$ | | | | | | | | |
| **8** | 8 | 3 | 6 | 5 | 26 | 15 | 35 | |
| **7** | 8 | 3 | 2 | 5 | 26 | 15 | | |
| **6** | 8 | 3 | 6 | 5 | 26 | | | |
| **5** | 8 | 3 | 3 | 5 | | | | |
| **4** | 8 | 3 | 6 | | | | | |
| **3** | 8 | 3 | | | | | | |
| **2** | 8 | | | | | | | |
| | **2** | **3** | **4** | **5** | **6** | **7** | **8** | λ |

| | | | t=10 | | | | | | |
|---|---|---|---|---|---|---|---|---|---|
| $b_0$ | | | | | | | | | |
| **9** | 18 | 28 | 6 | 42 | 58 | 22 | 13 | 44 | |
| **8** | 18 | 28 | 6 | 42 | 58 | 22 | 13 | | |
| | 18 | 28 | 6 | 6 | 58 | 22 | | | |
| **6** | 18 | 28 | 6 | 42 | 58 | | | | |
| **5** | 18 | 28 | 6 | 42 | | | | | |
| **4** | 18 | 28 | 6 | | | | | | |
| **3** | 18 | 28 | | | | | | | |
| **2** | 18 | | | | | | | | |
| | **2** | **3** | **4** | **5** | **6** | **7** | **8** | **9** | λ |

| | | | | | | | t=11 | | | |
|---|---|---|---|---|---|---|---|---|---|---|
| $b_0$ | | | | | | | | | | |
| **10** | 6 | 4 | 7 | 18 | 12 | 3 | 28 | 21 | 108 | |
| **9** | 3 | 8 | 7 | 2 | 12 | 6 | 28 | 21 | | |
| **8** | 6 | 2 | 7 | 18 | 12 | 3 | 28 | | | |
| **7** | 2 | 8 | 7 | 18 | 12 | 6 | | | | |
| **6** | 3 | 4 | 7 | 6 | 18 | | | | | |
| **5** | 6 | 8 | 7 | 18 | | | | | | |
| **4** | 6 | 2 | 7 | | | | | | | |
| **3** | 3 | 8 | | | | | | | | |
| **2** | 6 | | | | | | | | | |
| | **2** | **3** | **4** | **5** | **6** | **7** | **8** | **9** | **10** | λ |

| | | | | | t=12 | | | | | | |
|---|---|---|---|---|---|---|---|---|---|---|---|
| $b_0$ | | | | | | | | | | | |
| **11** | 11 | 12 | 23 | 29 | 35 | 41 | 12 | 53 | 48 | 63 | |
| **10** | 11 | 12 | 23 | 29 | 35 | 41 | 6 | 53 | 48 | | |
| **9** | 11 | 12 | 23 | 29 | 35 | 41 | 12 | 53 | | | |
| **8** | 11 | 12 | 23 | 29 | 35 | 41 | 12 | | | | |
| **7** | 11 | 12 | 23 | 29 | 35 | 41 | | | | | |
| **6** | 11 | 12 | 23 | 29 | 35 | | | | | | |
| **5** | 11 | 12 | 23 | 29 | | | | | | | |
| **4** | 11 | 12 | 23 | | | | | | | | |
| **3** | 11 | 12 | | | | | | | | | |
| **2** | 11 | | | | | | | | | | |
| | **2** | **3** | **4** | **5** | **6** | **7** | **8** | **9** | **10** | **11** | **λ** |

| | | | | | t=13 | | | | | | | |
|---|---|---|---|---|---|---|---|---|---|---|---|---|
| $b_0$ | | | | | | | | | | | | |
| **12** | 20 | 18 | 4 | 4 | 10 | 4 | 17 | 14 | 21 | 70 | 60 | |
| **11** | 20 | 18 | 4 | 16 | 2 | 12 | 17 | 14 | 21 | 70 | | |
| **10** | 4 | 18 | 4 | 8 | 10 | 3 | 17 | 14 | 21 | | | |
| **9** | 20 | 18 | 4 | 16 | 10 | 4 | 17 | 14 | | | | |
| **8** | 20 | 18 | 4 | 2 | 10 | 12 | 17 | | | | | |
| **7** | 20 | 18 | 4 | 16 | 10 | 12 | | | | | | |
| **6** | 20 | 18 | 4 | 8 | 10 | | | | | | | |
| **5** | 4 | 18 | 4 | 16 | | | | | | | | |
| **4** | 20 | 18 | 4 | | | | | | | | | |
| **3** | 20 | 18 | | | | | | | | | | |
| **2** | **20** | | | | | | | | | | | |
| | **2** | **3** | **4** | **5** | **6** | **7** | **8** | **9** | **10** | **11** | **12** | **λ** |

| | | | | | | | t=14 | | | | | | |
|---|---|---|---|---|---|---|---|---|---|---|---|---|---|
| $b_0$ | | | | | | | | | | | | | |
| **13** | 18 | 8 | 10 | 22 | 82 | 96 | 12 | 50 | 46 | 48 | 83 | 45 | |
| **12** | 6 | 8 | 10 | 22 | 82 | 96 | 12 | 50 | 46 | 16 | 83 | | |
| **11** | 18 | 8 | 2 | 22 | 82 | 96 | 12 | 50 | 46 | 48 | | | |
| **10** | 18 | 8 | 5 | 22 | 82 | 96 | 12 | 10 | 46 | | | | |
| **9** | 2 | 8 | 10 | 22 | 82 | 96 | 12 | 50 | | | | | |
| **8** | 18 | 8 | 10 | 22 | 82 | 96 | 12 | | | | | | |
| **7** | 18 | 8 | 10 | 22 | 82 | 96 | | | | | | | |
| **6** | 6 | 8 | 10 | 22 | 82 | | | | | | | | |
| **5** | 18 | 8 | 5 | 22 | | | | | | | | | |
| **4** | 18 | 8 | 10 | | | | | | | | | | |
| **3** | 6 | 8 | | | | | | | | | | | |
| **2** | 18 | | | | | | | | | | | | |
| | **2** | **3** | **4** | **5** | **6** | **7** | **8** | **9** | **10** | **11** | **12** | **13** | **λ** |

| | | | | | | | t=15 | | | | | | | |
|---|---|---|---|---|---|---|---|---|---|---|---|---|---|---|
| ***b*<sub></sub>** | | | | | | | | | | | | | | |
| **14** | 28 | 5 | 29 | 36 | 88 | 12 | 8 | 11 | 148 | 40 | 89 | 96 | 90 | |
| **13** | 28 | 10 | 29 | 36 | 88 | 2 | 8 | 11 | 148 | 40 | 89 | 96 | | |
| **12** | 28 | 5 | 29 | 36 | 88 | 12 | 8 | 11 | 148 | 40 | 89 | | | |
| **11** | 28 | 2 | 29 | 36 | 88 | 12 | 8 | 11 | 148 | 40 | | | | |
| **10** | 28 | 5 | 29 | 36 | 88 | 12 | 8 | 11 | 148 | | | | | |
| **9** | 28 | 10 | 29 | 36 | 88 | 12 | 8 | 11 | | | | | | |
| **8** | 28 | 5 | 29 | 36 | 88 | 12 | 8 | | | | | | | |
| **7** | 28 | 10 | 29 | 36 | 88 | 12 | | | | | | | | |
| **6** | 28 | 5 | 29 | 36 | 88 | | | | | | | | | |
| **5** | 28 | 10 | 29 | 36 | | | | | | | | | | |
| **4** | 28 | 5 | 29 | | | | | | | | | | | |
| **3** | 28 | 10 | | | | | | | | | | | | |
| **2** | 28 | | | | | | | | | | | | | |
| | **2** | **3** | **4** | **5** | **6** | **7** | **8** | **9** | **10** | **11** | **12** | **13** | **14** | **λ** |

| | | | | | | | t=16 | | | | | | | | |
|---|---|---|---|---|---|---|---|---|---|---|---|---|---|---|---|
| $b_0$ | | | | | | | | | | | | | | | |
| **15** | 5 | 23 | 3 | 39 | 9 | 9 | 7 | 3 | 13 | 3 | 95 | 11 | 37 | 119 | |
| **14** | 5 | 23 | 3 | 39 | 9 | 9 | 7 | 15 | 13 | 5 | 95 | 33 | 37 | | |
| **13** | 5 | 23 | 3 | 39 | 9 | 9 | 7 | 5 | 13 | 15 | 95 | 33 | | | |
| **12** | 5 | 23 | 3 | 39 | 9 | 9 | 7 | 15 | 13 | 15 | 95 | | | | |
| **11** | 5 | 23 | 3 | 39 | 9 | | 7 | 3 | 13 | 15 | | | | | |
| **10** | 5 | 23 | 3 | 39 | 9 | 9 | 7 | 15 | 13 | | | | | | |
| **9** | 5 | 23 | 3 | 39 | 9 | 9 | 7 | 15 | | | | | | | |
| **8** | 5 | 23 | 3 | 39 | 9 | 9 | 7 | | | | | | | | |
| **7** | 5 | 23 | 3 | 39 | 9 | 9 | | | | | | | | | |
| **6** | 5 | 23 | 3 | 39 | 9 | | | | | | | | | | |
| **5** | 5 | 23 | 3 | 39 | | | | | | | | | | | |
| **4** | 5 | 23 | 3 | | | | | | | | | | | | |
| **3** | 5 | 23 | | | | | | | | | | | | | |
| **2** | 5 | | | | | | | | | | | | | | |
| | **2** | **3** | **4** | **5** | **6** | **7** | **8** | **9** | **10** | **11** | **12** | **13** | **14** | **15** | **λ** |

**t=17**

| $b_0$ | | | | | | | | | | | | | | | | |
|---|---|---|---|---|---|---|---|---|---|---|---|---|---|---|---|---|
| **16** | 10 | 20 | 33 | 6 | 10 | 29 | 12 | 9 | 78 | 30 | 12 | 20 | 26 | 63 | 135 | |
| **15** | 10 | 4 | 33 | 6 | 10 | 29 | 2 | 9 | 78 | 30 | 12 | 10 | 26 | 63 | | |
| **14** | 10 | 20 | 33 | 2 | 10 | 29 | 12 | 9 | 78 | 30 | 4 | 20 | 26 | | | |
| **13** | 10 | 20 | 33 | 6 | 10 | 29 | 12 | 9 | 6 | 30 | 12 | 20 | | | | |
| **12** | 10 | 20 | 33 | 6 | 10 | 29 | 4 | 9 | 78 | 30 | 12 | | | | | |
| **11** | 2 | 20 | 33 | 6 | 10 | 29 | 12 | 9 | 78 | 30 | | | | | | |
| **10** | 10 | 4 | 33 | 6 | 10 | 29 | 6 | 9 | 78 | | | | | | | |
| **9** | 10 | 20 | 33 | 6 | 10 | 29 | 4 | 9 | | | | | | | | |
| **8** | 10 | 20 | 33 | 6 | 10 | 29 | 12 | | | | | | | | | |
| **7** | 10 | 20 | 33 | 2 | 10 | 29 | | | | | | | | | | |
| **6** | 10 | 20 | 33 | 6 | 10 | | | | | | | | | | | |
| **5** | 10 | 4 | 33 | 6 | | | | | | | | | | | | |
| **4** | 10 | 20 | 33 | | | | | | | | | | | | | |
| **3** | 10 | 20 | | | | | | | | | | | | | | |
| **2** | 10 | | | | | | | | | | | | | | | |
| | **2** | **3** | **4** | **5** | **6** | **7** | **8** | **9** | **10** | **11** | **12** | **13** | **14** | **15** | **16** | **λ** |

**t=18**

| $b_0$ | | | | | | | | | | | | | | | | | |
|---|---|---|---|---|---|---|---|---|---|---|---|---|---|---|---|---|---|
| **17** | 12 | 52 | 35 | 44 | 106 | 20 | 20 | 33 | 178 | 196 | 84 | 116 | 250 | 262 | 15 | 60 | |
| **16** | 12 | 52 | 35 | 44 | 106 | 20 | 20 | 33 | 178 | 196 | 84 | 116 | 250 | 262 | 15 | | |
| **15** | 3 | 52 | 35 | 44 | 106 | 4 | 20 | 33 | 178 | 196 | 42 | 116 | 250 | 262 | | | |
| **14** | 4 | 52 | 35 | 44 | 106 | 20 | 20 | 11 | 178 | 196 | 84 | 116 | 250 | | | | |
| **13** | 12 | 52 | 35 | 44 | 106 | 20 | 10 | 33 | 178 | 196 | 84 | 116 | | | | | |
| **12** | 12 | 52 | 35 | 44 | 106 | 20 | 20 | 33 | 178 | 196 | 84 | | | | | | |
| **11** | 12 | 52 | 35 | 44 | 106 | 20 | 4 | 33 | 178 | 196 | | | | | | | |
| **10** | 3 | 52 | 35 | 44 | 106 | 4 | 20 | 33 | 178 | | | | | | | | |
| **9** | 12 | 52 | 35 | 44 | 106 | 20 | 20 | 33 | | | | | | | | | |
| **8** | 12 | 52 | 35 | 44 | 106 | 20 | 20 | | | | | | | | | | |
| **7** | 4 | 52 | 35 | 44 | 106 | 20 | | | | | | | | | | | |
| **6** | 12 | 52 | 35 | 44 | 106 | | | | | | | | | | | | |
| **5** | 3 | 52 | 35 | 44 | | | | | | | | | | | | | |
| **4** | 12 | 52 | 35 | | | | | | | | | | | | | | |
| **3** | 12 | 52 | | | | | | | | | | | | | | | |
| **2** | 12 | | | | | | | | | | | | | | | | |
| | **2** | **3** | **4** | **5** | **6** | **7** | **8** | **9** | **10** | **11** | **12** | **13** | **14** | **15** | **16** | **17** | **λ** |

**t=19**

| $b_0$ | | | | | | | | | | | | | | | | | |
|---|---|---|---|---|---|---|---|---|---|---|---|---|---|---|---|---|---|
| **18** | 36 | 6 | 10 | 46 | 112 | 10 | 5 | 8 | 6 | 12 | 113 | 40 | 52 | 35 | 25 | 66 | 30 | |
| **17** | 36 | 6 | 10 | 46 | 112 | 10 | 5 | 2 | 6 | 12 | 113 | 40 | 52 | 70 | 25 | 66 | | |
| **16** | 36 | 6 | 10 | 46 | 112 | 10 | 5 | 8 | 6 | 12 | 113 | 40 | 52 | 35 | 25 | | | |
| **15** | 36 | 6 | 2 | 46 | 112 | 10 | 5 | 8 | 6 | 12 | 113 | 40 | 52 | 70 | | | | |
| **14** | 36 | 2 | 10 | 46 | 112 | 10 | 5 | 8 | 3 | 12 | 113 | 40 | 52 | | | | | |
| **13** | 36 | 6 | 10 | 46 | 112 | 10 | 5 | 8 | 6 | 4 | 113 | 40 | | | | | | |
| **12** | 36 | 6 | 10 | 46 | 112 | 10 | 5 | 8 | 6 | 12 | 113 | | | | | | | |
| **11** | 36 | 6 | 10 | 46 | 112 | 4 | 5 | 8 | 6 | 12 | | | | | | | | |
| **10** | 36 | 6 | 2 | 46 | 112 | 10 | 5 | 8 | 6 | | | | | | | | | |
| **9** | 36 | 6 | 10 | 46 | 112 | 10 | 5 | 8 | | | | | | | | | | |
| **8** | 36 | 6 | 10 | 46 | 112 | 10 | 5 | | | | | | | | | | | |
| **7** | 36 | 2 | 10 | 46 | 112 | 10 | | | | | | | | | | | | |
| **6** | 36 | 6 | 10 | 46 | 112 | | | | | | | | | | | | | |
| **5** | 36 | 6 | 2 | 46 | | | | | | | | | | | | | | |
| **4** | 36 | 6 | 10 | | | | | | | | | | | | | | | |
| **3** | 36 | 6 | | | | | | | | | | | | | | | | |
| **2** | 36 | | | | | | | | | | | | | | | | | |
| | **2** | **3** | **4** | **5** | **6** | **7** | **8** | **9** | **10** | **11** | **12** | **13** | **14** | **15** | **16** | **17** | **18** | **λ** |

**t=20**

| $b_0$ | | | | | | | | | | | | | | | | | | |
|---|---|---|---|---|---|---|---|---|---|---|---|---|---|---|---|---|---|---|
| **19** | 12 | 29 | 39 | 30 | 16 | 69 | 52 | 89 | 99 | 72 | 119 | 36 | 30 | 132 | 35 | 112 | 179 | 189 | |
| **18** | 12 | 29 | 39 | 5 | 16 | 69 | 52 | 89 | 99 | 72 | 119 | 36 | 15 | 132 | 35 | 112 | 179 | | |
| **17** | 12 | 29 | 39 | 30 | 2 | 69 | 52 | 89 | 99 | 72 | 119 | 36 | 30 | 132 | 35 | 112 | | | |
| **16** | 12 | 29 | 39 | 30 | 16 | 69 | 52 | 89 | 99 | 72 | 119 | 36 | 30 | 132 | 35 | | | | |
| **15** | 12 | 29 | 39 | 10 | 16 | 69 | 52 | 89 | 99 | 72 | 119 | 36 | 30 | 132 | | | | | |
| **14** | 12 | 29 | 39 | 30 | 16 | 69 | 52 | 89 | 99 | 72 | 119 | 36 | 30 | | | | | | |
| **13** | 12 | 29 | 39 | 30 | 16 | 69 | 52 | 89 | 99 | 72 | 119 | 36 | | | | | | | |
| **12** | 12 | 29 | 39 | 10 | 16 | 69 | 52 | 89 | 99 | 72 | 119 | | | | | | | | |
| **11** | 12 | 29 | 39 | 6 | 16 | 69 | 52 | 89 | 99 | 72 | | | | | | | | | |
| **10** | 12 | 29 | 39 | 30 | 16 | 69 | 52 | 89 | 99 | | | | | | | | | | |
| **9** | 12 | 29 | 39 | 5 | 16 | 69 | 52 | 89 | | | | | | | | | | | |
| **8** | 12 | 29 | 39 | 30 | 16 | 69 | 52 | | | | | | | | | | | | |
| **7** | 12 | 29 | 39 | 30 | 16 | 69 | | | | | | | | | | | | | |
| **6** | 12 | 29 | 39 | 10 | 16 | | | | | | | | | | | | | | |
| **5** | 12 | 29 | 39 | 30 | | | | | | | | | | | | | | | |
| **4** | 12 | 29 | 39 | | | | | | | | | | | | | | | | |
| **3** | 12 | 29 | | | | | | | | | | | | | | | | | |
| **2** | 12 | | | | | | | | | | | | | | | | | | |
| | **2** | **3** | **4** | **5** | **6** | **7** | **8** | **9** | **10** | **11** | **12** | **13** | **14** | **15** | **16** | **17** | **18** | **19** | **λ** |

# APENDIX I. Transformation of the congruence (II.5) into (II.6).

To transform congruence **(II.5) into (II.6**) we use **General Divisibility Criteria (GDC)** derived in ref.[3]. The **GDC** of any number **A** (in the form **(II.1))** by divider of the form **p=(wt – u)**, is expressed by

$$\mathbf{GDC}(\boldsymbol{B}) = \sum_{k=0}^{N-1} \boldsymbol{u}^{k} \boldsymbol{w}^{N-k}\, \boldsymbol{b}_{k} \qquad \textbf{(A.1)}$$

where $b_k$ are digits of **B** in the **t**-base number system. Note, that parameters **w** and **u** can be written for arbitrary number**.** In agreement with meaning of the divisibility criteria, the numbers **B** and **GDC(B)** must both be divisible or not divisible by **p**. Despite the universality **GDC** it is not always more simple than number **B**. This depends on the digits of number **B** as well as on divider parameters **w** and **u**. In the case of Eq.**(II.5)** the all **N** digits of denominator are the same and equal to $b_k = b_0(t-1)$ . The parameters of the divider **p** are also simple: **w=** $\lambda$ , **u=1**. Substituting these numbers into **(A.1)** we find

$$\mathbf{GDC}(P_{\lambda}^{N}) = b_0\,(\mathrm{t}-1) \sum_{k=0}^{N-1} \lambda^{k} = b_0\,(t-1)\,\frac{(\lambda^{N}-1)}{(\lambda-1)} \qquad \textbf{(A.2)}$$

The divisibility of $\mathbf{GDC}(P_{\lambda}^{N})$ by **p** is expressed through the congruence

$$b_0 \frac{(t-1)}{\lambda-1}(\lambda^{N}-1) \equiv 0 \qquad (mod\,(\lambda\, t-1)) \qquad \textbf{(A.3)}$$

We will show that multiplier **β =(t-1)/(λ-1)** is co-prime to module **p.** In this case **β** can be removed from **(A.2)** without altering the module. First, we note that **β** is smaller than **p** and therefore the module cannot be contained in **β**. This follows from the inequalities **β** $< (t-1) < (\lambda t - 1)$**.** Further, presenting module in the form $p = \lambda * (t-1) + (\lambda - 1)$ we note that if **p** and **(t-1)** are both multiples of some number **δ**, than $(\lambda - 1)$ is also multiple of **δ**. Therefore, the **δ** is contracted in numerator and denominator of **β.** Thus **β** is co-prime to **p** and can be removed from the congruence. This leads to the congruence **(II.6).**

## Conclusion

The formula for parasitic numbers was derived based on their definition. Its application is demonstrated on several number systems. In conclusion we'd like to say several worlds about congruence **(II.6)** or **(II.7** solution. They are related to the problem of power residues elaborated in details by Gauss [4]. It acquires the simplest form when module p is prime number. This occurs in many cases that were solved above. Nevertheless, in our specific problem the numerical solution is very simple and in many cases, so simple that even a solution by hand is available, using properties of residues.

## References

[1] John Tierney "Prize for Dyson Puzzle", New York Times, April 13, 2009.

[2] See for instance Wikipedia.

[3] Anatoly A. Grinberg and Serge Luryi, "General Divisibility Criteria," available (since Jan 21, 2014) online at arXiv.org : 1401.5486.

[4] Carl Friedrich Gauss « Disquisitiones arithmetic».